\documentclass[12pt]{article}

\usepackage[cp1250]{inputenc}
\usepackage{amsmath}
\usepackage{amssymb}
\usepackage[a4paper,left=3cm,textwidth=15cm]{geometry}

\newcommand{\p}{\mathbb{P}}
\newcommand{\N}{\mathbb{N}}
\newcommand{\E}{\mathbb{E}}
\newcommand{\Var}{\mathrm{Var}}
\newcommand{\R}{\mathbb{R}}

\newcommand{\proof}{\noindent \textbf{Proof. \quad}}

\newcommand{\remark}{\paragraph{Remark}}
\newcommand{\qed}{\begin{flushright}\nopagebreak[4]$\square$\end{flushright}}
\newcommand{\ind}{\mathbf{1}}
\newcommand{\Ent}{\mathrm{Ent}\,}

\newtheorem{theorem}{Theorem}
\newtheorem{lemma}{Lemma}
\newtheorem{cor}{Corollary}
\newtheorem{defi}{Definition}

\newtheorem{prop}{Proposition}

\author{Rados\l aw Adamczak}
\title{Logarithmic Sobolev inequalities and concentration of measure for convex functions and polynomial chaoses.}
\date{}
\begin{document}
\maketitle
\begin{abstract}
We prove logarithmic Sobolev inequalities and concentration
results for convex functions and a class of product random
vectors. The results are used to derive tail and moment
inequalities for chaos variables (in the spirit of Talagrand and
Arcones, Gin\'{e}). We also show that the same proof may be used
for chaoses generated by log-concave random variables, recovering
results by \L ochowski and present an application to exponential
integrability of Rademacher chaos.
\end{abstract}

\section{Introduction}
The paper is concerned with concentration properties of random
vectors. We start with the following
\begin{defi}
A real random variable $\xi$ is said to have \emph{the
concentration property} of order $\alpha > 0$ with constants $K,C$
if there exists $a \in \R$ such that for all $t \ge 0$,
\begin{equation}
\p(|\xi - a| \ge t) \le Ce^{-t^\alpha/K}.
\end{equation}
\end{defi}

It is easy to see that the concentration property implies that
$\xi$ has a finite moment and there exist constants $C',K'$,
depending on $\alpha,C,K$ only, such that for all $t \ge 0$,
\begin{equation} \label{alfa_mean_concentration}
\p(|\xi - \E \xi| \ge t) \le C'e^{-t^\alpha/K'}.
\end{equation}

Moreover, by the Chebyshev inequality, the condition
(\ref{alfa_mean_concentration}) is equivalent to the following
moment estimates, valid for all $p \ge 1$:
\begin{equation} \label{alfa_moment_concentration}
\|\xi - \E \xi\|_p \le K''p^{1/\alpha}.
\end{equation}
More precisely, if (\ref{alfa_mean_concentration}) holds then so
does (\ref{alfa_moment_concentration}) with $K''$ depending only
on $\alpha,C',K'$, whereas (\ref{alfa_moment_concentration})
implies (\ref{alfa_mean_concentration}) with $C',K'$ depending
only on $K'',\alpha$.

In what follows we restrict our attention to random variables of
the form $\xi = f(X)$, where $X$ is a random vector in $\R^n$ and
the function $f$ belongs to $\mathcal{F}$,  a specified class of
real, Borel measurable functions on $\R^n$ (e.g. $1$-Lipschitz
functions or $1$-Lipschitz (homogeneous) convex functions).

\begin{defi}
We say that a random vector $X$ in $\R^n$ has \emph{the
concentration property} of order $\alpha$ with constants $C,K$
with respect to a class $\mathcal{F}$ of real, Borel measurable
functions on $\R^n$ if for every $f\in \mathcal{F}$ the random
variable $f\circ X$ has the concentration property of order
$\alpha$ with constants $C,K$.
\end{defi}

The above definition seems justified as there are quite a few
examples of pairs $(X,\mathcal{F})$ satisfying it. For instance,
it is by now classical that if $X$ is a standard Gaussian random
vector in $\R^n$ then it has the concentration property of order 2
with constants 1,2 with respect to the class of $1$-Lipschitz
functions. Also random vectors in $\R^n$ with independent
uniformly bounded components have the concentration property of
order 2 with constants independent of the dimension with
$\mathcal{F}$ being the class of $1$-Lipschitz convex functions
\cite{T3}. The latter example can be extended to arbitrary random
vectors with bounded support, but the constants will then also
depend on the mixing coefficients associated with the random
vector \cite{S}.

We now briefly describe one of the most efficient tools for
proving the concentration property (especially for product
distributions), which has been developed over the past several
years, namely the entropy method.
\begin{defi}
Let $\xi$ be a nonnegative random variable and $\Phi \colon \R_+
\to \R$ a convex function such that $\Phi''>0$ and $1/\Phi''$ is
concave. Define \emph{the $\Phi$-entropy} of $\xi$ by the formula
\begin{equation}
\Ent_\Phi \xi = \E\Phi(\xi) - \Phi(\E \xi).
\end{equation}
\end{defi}

The most important examples are $\Phi(x) = x^2$ and $\Phi(x) =
x\log x$. In these cases $\Ent_\Phi$ becomes respectively the
variance and the usual entropy of a random variable (which will be
denoted simply by $\Ent$). The notion is important from the
concentration of measure point of view since we have
\begin{theorem}[Herbst argument, see
\cite{L1},\cite{L2}]\label{log_sobolev_theorem} Let $X \in \R^n$
be a random variable and $\mathcal{F}$ a class of functions such
that $\lambda f \in \mathcal{F}$ for all $f \in \mathcal{F}$ and
$\lambda \ge 0$. Assume furthermore that for all $f \in
\mathcal{F}$,
\begin{equation}\label{log_sobolev_inequality}
\Ent e^{f(X)} \le C\E|\nabla f(X)|^2e^{f(X)},
\end{equation}
and the right-hand side is finite. Then for all $f \in
\mathcal{F}$ with $|\nabla f| \le 1$ and $t \ge 0$,
\begin{displaymath}
\p(f(X) \ge \E f(X) + t) \le e^{-t^2/4C}.
\end{displaymath}
\end{theorem}

A crucial property of $\Ent_\Phi$ is \emph{the tensorization},
which is described in the following
\begin{theorem}[see \cite{BBLM},\cite{LO}]\label{tensorization_theorem} Consider a
product probability space $(\Omega, \mu)$, where $\Omega =
\bigotimes_{i=1}^n \Omega_i$ and $\mu = \otimes_{i=1}^n\mu_i$.
Then for every nonnegative random variable $\xi$ we have
\begin{displaymath}
\Ent_{\Phi}\xi \le \sum_{i = 1}^{n} \E\; \Ent_{\Phi,\mu_i}\xi,
\end{displaymath}
where $\Ent_{\Phi,\mu_i}\xi$ denotes the value of the functional
$\Ent_{\Phi}$ at the function $\xi$, considered as a function on
$\Omega_i$, with the other coordinates fixed. \label{tensor}
\end{theorem}

Thus if a random vector $X \in \R^n$ satisfies the inequality
(\ref{log_sobolev_inequality}) for all $f \in \mathcal{F}$, then
so does the random vector $X_1\oplus\ldots\oplus X_d \in
(\R^n)^d$, where $X_i$ are independent  copies of $X$, for all
functions $f \colon (\R^n)^d \to \R$ such that
$f(x_1,\ldots,x_{i-1},\cdot,x_{i+1},\ldots,x_d) \in \mathcal{F}$
for all $i$ and $x_1,\ldots,x_n \in \R^n$, which can be used to
obtain concentration inequalities. This method has led to
concentration results for $1$-Lipschitz functions of standard
Gaussian vectors and $1$-Lipschitz convex functions of uniformly
bounded product distributions (see \cite{L2}, chapter 5). In a
slightly different setting it was also used to obtain
concentration results for more general functions of independent
random variables and also to some general moment inequalities for
such functions \cite{BBLM}. We also mention that inequalities in
the spirit of (\ref{log_sobolev_inequality}) with the left-hand
side replaced by $\Var f(X)$ (the so called Poincar\'{e}
inequalities) yield concentration property of order 1. There are
also similar Lata\l a-Oleszkiewicz inequalities which imply
concentration of order $\alpha \in (1,2)$ \cite{LO}.

In this article we will present two results concerning
concentration. First, in Section \ref{Sec_log_Sob} we obtain some
sufficient conditions for a real random variable to satisfy the
logarithmic Sobolev inequality (\ref{log_sobolev_inequality}) for
convex functions, which yields some subgaussian deviation
inequalities. Then in Section \ref{Sec_tensor_product} we will
show that the concentration property of a random vector $X$ with
respect to the class of semi-norms can be tensorized to obtain
concentration inequalities for $X_1\otimes\ldots\otimes X_d$
(where $X_i$'s are independent copies of $X$), which gives some
new and helps to recover known inequalities for polynomial
chaoses. Finally, in the last section we present an application of
these inequalities, by presenting a new proof of exponential
integrability for Rademacher chaos process.

\section{Logarithmic Sobolev inequalities and concentration of measure for convex
functions \label{Sec_log_Sob}}

\begin{defi}\label{class_m}
For $m > 0$ and $\sigma \ge 0$ let $\mathcal{M}(m,\sigma^2)$
denote the class of probability distributions $\mu$ on $\R$ for
which
\begin{displaymath}
\upsilon^+(A) \le \sigma^2\mu(A)
\end{displaymath}
for all sets $A$ of the form $A = [x,\infty)$, $x \ge m$ and
\begin{displaymath}
 \upsilon^-(A) \le \sigma^2 \mu(A)
\end{displaymath}
for all sets $A$ of the form $A = (-\infty,-x]$,  $x \ge m$, where
$\upsilon^+$ is the measure on $[m,\infty)$ with density $g(x) =
x\mu([x,\infty))$ and $\upsilon^-$ is the measure on
$(-\infty,-m]$ with density $g(x) = -x\mu((-\infty,x])$.
\end{defi}

\begin{prop}\label{equivalence}
Let $\mu$ be a probability distribution on $\R$. Then the
following conditions are equivalent:
\begin{itemize}
\item[(i)] $\mu \in \mathcal{M}(m,\sigma^2)$ for some $m,\sigma$,

\item[(ii)]
\begin{eqnarray*}
\mu([x+ C/x,\infty)) &\le& \alpha\mu([x,\infty)) \\
\mu((-\infty,-x - C/x]) &\le& \alpha\mu((-\infty,-x])
\end{eqnarray*}
for some $C > 0$, $\alpha < 1$ and all $x \ge m$,
\end{itemize}
where the constants in (ii) depend only on the constants in (i)
and vice versa. For instance if (i) holds, we can take $C =
2\sigma^2$ and $\alpha = 1/2$.
\end{prop}

\proof Assume (i) holds. Then for $x \ge m$ we have
\begin{displaymath}
\sigma^2 \mu([x,\infty)] \ge \int_{x}^{ x + 2\sigma^2/x}y
\mu([y,\infty))dy \ge
x\frac{2\sigma^2}{x}\mu([x+2\sigma^2/x,\infty)) =
2\sigma^2\mu([x+2\sigma^2/x,\infty)),
\end{displaymath}
which clearly implies the first inequality of (ii). The second
inequality follows similarly.

Suppose now that (ii) is satisfied and for $x \ge m$ define the
sequence $a_0 = x, a_{n+1} = a_n + C/a_n$. Then it is easy to see
that $a_n \to \infty$ and therefore
\begin{eqnarray*}
\int_x^\infty t\mu([t,\infty))dt \le \sum_{n=0}^\infty
a_{n+1}\mu([a_n,\infty))(a_{n+1} - a_n) \le K_1 \sum_{n=0}^\infty
\alpha^n \mu([a_0,\infty)) \le K_2 \mu([x,\infty)).
\end{eqnarray*}
We can proceed analogously to obtain the condition on the left
tail. \qed

\remark It is also worth noting that for a real random variable
$X$, the condition
\mbox{$\mathcal{L}(X)\in\mathcal{M}(m,\sigma^2)$} is equivalent to
$\E X^2\ind_{\{X \ge t\}} \le (t^2 + 2\sigma^2)\p(X \ge t)$ and
$\E X^2\ind_{\{X \le -t\}} \le (t^2 + 2\sigma^2)\p(X \le -t)$ for
all $t \ge m$.

\paragraph{Example} Of course all measures with bounded support
belong to $\mathcal{M}(m,0)$ for some $m$. Other examples of
measures from $\mathcal{M}(m,\sigma^2)$ are absolutely continuous
distributions $\mu$ satisfying the inequalities
\begin{displaymath}
\frac{d}{dt}\log \mu([t,\infty)) \le - \frac{t}{\sigma^2}, \quad
\frac{d}{dt}\log \mu((-\infty,-t]) \le -\frac{t}{\sigma^2}
\end{displaymath}
for $t \ge m$. In particular if $\mu$ has density of the form
$e^{-V(x)}$ with $V'(x) \ge x/\sigma^2$ and $V'(-x) \le
-x/\sigma^2$ then $\mu \in \mathcal{M}(1,\sigma^2)$.

\paragraph{}Now we are ready to state the main result of this section.

\begin{theorem}\label{entropy_estimates}
Let $X_1,\ldots,X_n$ be independent random variables such that
$\mathcal{L}(X_i) \in \mathcal{M}(m,\sigma^2)$ and let
$\varphi\colon \R^n \to \R$ be a smooth convex function. Then
\begin{displaymath}
\Ent e^{\varphi(X_1,\ldots,X_n)} \le C(m,\sigma^2)\E
e^{\varphi(X_1,\ldots,X_n)}|\nabla \varphi(X_1,\ldots,X_n)|^2.
\end{displaymath}
Hence for every $1$-Lipschitz convex function $\varphi\colon \R^n
\to \R$ and all $t \ge 0$,
\begin{displaymath}
\p(\varphi(X_1,\ldots,X_n) \ge  \E \varphi(X_1,\ldots,X_n) + t)
\le e^{-\frac{t^2}{4C(m,\sigma^2)}}.
\end{displaymath}
\end{theorem}
Before we proceed to  the proof of Theorem
\ref{entropy_estimates}, we will need a few lemmas.
\begin{lemma}\label{double_monotone_function}
Let $\mu \in \mathcal{M}(m,\sigma^2)$. Then for all functions $f
\colon \R \to \R_+$ which are non-increasing for $x \le x_0$ and
non-decreasing for $x \ge x_0$, we have
\begin{displaymath}
\int_{\tilde{m}}^\infty f(x)x \mu([x,\infty))dx \le 2\sigma^2
\int_\R f(x)d\mu(x),
\end{displaymath}
where $\tilde{m} = m\vee(\sqrt{2}\sigma) +
2\sigma^2/(m\vee(\sqrt{2}\sigma))$.
\end{lemma}
\proof First notice that by standard approximation arguments the
inequalities of Definition \ref{class_m} are also satisfied for
sets $A = (x,\infty)$, $x \ge m$. We have

\begin{eqnarray} \label{layer_cake}
&&\int_{\tilde{m}}^\infty f(x)x\mu([x,\infty))dx =
\int_{\tilde{m}}^\infty\int_0^\infty \ind_{\{s \le
f(x)\}}x\mu([x,\infty))dsdx \nonumber\\
&=& \int_0^\infty\upsilon^+(\{x \ge \tilde{m} \colon f(x) \ge
s\})ds.
\end{eqnarray}
The set $A = \{x \ge \tilde{m}\colon f(x) \ge s\}$ is either a
semi-line contained in $[m,\infty)$ or a disjoint union of such a
semi-line and an interval $I$ with the left end equal $\tilde{m}$.
In the former case we have $\upsilon^+(A) \le \sigma^2 \mu(A) \le
\sigma^2\mu(\{x\in \R\colon f(x) \ge s\})$.

Now consider the latter case. Denote the right end of $I$ by $t$.
Let $a = m\vee(\sqrt{2}\sigma)$. As $t \ge \tilde{m} = a +
2\sigma^2/a \ge a$, we obtain
\begin{eqnarray*}
\upsilon^+(A) &\le& \upsilon^+([\tilde{m},\infty)) \le
\sigma^2\mu([a,\infty)) \le 2\sigma^2\mu([a,a + 2\sigma^2/a))\\
&\le& 2\sigma^2 \mu([a,t)) \le 2\sigma^2 \mu(\{x \le t \colon f(x)
\ge s\}),
\end{eqnarray*}
where the second inequality follows from the assumption $\mu \in
\mathcal{M}(m,\sigma^2)$, the third from Proposition
\ref{equivalence} and the last one from the observation that $x_0
\ge t$ and thus $f$ is nonincreasing for $x \le t$. Now we can
write
\begin{displaymath}
2\sigma^2\int_\R f(x)d\mu(x) = 2\sigma^2\int_\R \int_0^\infty
\ind_{\{s \le f(x)\}}ds d\mu(x) = 2\sigma^2\int_0^\infty \mu(\{x
\in \R \colon f(x) \ge s\}) ds,
\end{displaymath}
which together with (\ref{layer_cake}) allows us to complete the
proof.
 \qed

\begin{lemma}\label{subgaussian}
If $X$ is a random variable such that $\mathcal{L}(X) \in
\mathcal{M}(m,\sigma^2)$ then \\$\p(|X| \ge t)\le
C_1(m,\sigma^2)e^{-t^2/C_2(m,\sigma^2)}$ for all $t \ge 0$.
\end{lemma}
\proof Obviously it is sufficient to prove the inequality for
$t\ge 4m$. Define $g(x) = \int_x^\infty y\p(X \ge y)dy$. Then for
$x \ge m$,
\begin{displaymath}
x g(x) \le x \sigma^2 \p(X \ge x),
\end{displaymath}
and thus
\begin{displaymath}
\int_z^\infty xg(x)dx \le \sigma^2 g(z), \; z \ge m.
\end{displaymath}
Let $f(z) = \int_z^\infty xg(x)dx$. Since the function $x \mapsto
xg(x)$ is continuous, we can rewrite the above inequality as
\begin{displaymath}
f'(z) \le -\frac{1}{\sigma^2}z f(z),\; z \ge m,
\end{displaymath}
which gives $f(z) \le C\exp(-z^2/2\sigma^2)$ with $C$ depending
only on $m$ and $\sigma^2$. Now, as $g$ is nonincreasing, for $z
\ge m$ we have
\begin{displaymath}
g(2z) \le \frac{1}{z^2}f(z) \le \frac{C}{z^2}e^{-z^2/2\sigma^2}
\end{displaymath}
and similarly
\begin{displaymath}
\p(X \ge 4x) \le \frac{g(2x)}{4x^2}\le \frac{C}{4
x^4}e^{-x^2/2\sigma^2}
\end{displaymath}
 for $x\ge m$.

The lower  tail can be dealt with analogously. \qed

\begin{lemma}
\label{positive_axes_lemma} Let $\varphi \colon \R \to \R$ be a
smooth convex Lipschitz function and $X$ a random variable with
$\mathcal{L}(X) \in \mathcal{M}(m,\sigma^2)$. Then there exists a
constant $C(m,\sigma^2)$ such that
\begin{equation}
\label{positive_axes_inequality}
 \int_0^\infty\int_0^\infty
\varphi'(x)\varphi'(y)e^{\varphi(y)}\p(X \le x\wedge y)\p(X \ge x
\vee y)dx dy \le C(m,\sigma^2) \E\varphi'(X)^2e^{\varphi(X)}.
\end{equation}
\end{lemma}

\proof Let us first notice that the left-hand side of
(\ref{positive_axes_inequality}) is equal to
\begin{displaymath}
 \int_0^\infty\int_0^y
\varphi'(x)\varphi'(y)(e^{\varphi(x)} + e^{\varphi(y)})\p(X \le
x)\p(X \ge y)dx dy.
\end{displaymath}

Since $\varphi$ is convex, there exists a point $x_0$ (possibly
$0$ or infinity) at which $\varphi$ attains its minimum on
$[0,\infty]$. Moreover $\varphi$ is non-increasing on $(0,x_0)$
and non-decreasing on $(x_0,\infty)$. Therefore for $x,y \in
(0,x_0)$ with $x \le y$ one has
\begin{displaymath}
\varphi'(x)\varphi'(y)(e^{\varphi(x)} +e^{\varphi(y)}) \le
2\varphi'(x)^2e^{\varphi(x)}.
\end{displaymath}
Thus for $\tilde{m}$ being the constant defined in Lemma
\ref{double_monotone_function} we have
\begin{eqnarray}
&&\int_0^{x_0\wedge \tilde{m}}
\int_0^y\varphi'(x)\varphi'(y)(e^{\varphi(x)} +e^{\varphi(y)})\p(X
\le x)\p(X \ge y) dxdy \nonumber\\
&\le&
2\int_0^{x_0\wedge \tilde{m}}\varphi'(x)^2 e^{\varphi(x)}\p(X \le x) \int_x^{x_0\wedge \tilde{m}} \p(X \ge y)dy dx \nonumber \\
&\le& 2\tilde{m} \E \int_0^{x_0\wedge \tilde{m}} \varphi'(x)^2
e^{\varphi(x)}\ind_{\{X \le x\}} \le
2\tilde{m}^2\E\varphi'(X)^2e^{\varphi(X)}, \label{positive_axes_1}
\end{eqnarray}
where the last inequality follows from the fact that if $x_0 > 0$,
then $\varphi'(x)^2 e^{\varphi(x)} \le
\varphi'(y)^2e^{\varphi(y)}$ for $y \le x \le x_0$.

On the other hand, for $x_0< x < y$ we have
$\varphi'(x)\varphi'(y)(e^{\varphi(x)} + e^{\varphi(y)}) \le
2\varphi'(y)^2e^{\varphi(y)}$. Obviously this is also the case if
$x < x_0 < y$, so
\begin{eqnarray}
&& \int_{x_0\vee
\tilde{m}}^\infty\int_0^y\varphi'(x)\varphi'(y)(e^{\varphi(x)}
+e^{\varphi(y)})\p(X \le x)\p(X \ge y) dxdy\nonumber\\
&\le& 2\int_{x_0\vee \tilde{m}}^\infty\varphi'(y)^2
e^{\varphi(y)}y\p(X \ge y)dy \le 4\sigma^2
\E\varphi'(X)^2e^{\varphi(X)},\nonumber\\
\label{positive_axes_2}
\end{eqnarray}
by Lemma \ref{double_monotone_function}, since $\mathcal{L}(X) \in
\mathcal{M}(m,\sigma^2)$.

So it remains to estimate the integral over the interval
$(x_0\wedge \tilde{m},x_0\vee \tilde{m})$. Let us consider two
cases.

\paragraph{(i)} $x_0 < \tilde{m}$, then
\begin{eqnarray}
&&\int_{x_0}^{\tilde{m}}\int_0^y
\varphi'(x)\varphi'(y)(e^{\varphi(x)} +
e^{\varphi(y)})\p(X \le x)\p(X \ge y)dx dy \nonumber \\
&\le& 2\int_{x_0}^{\tilde{m}}
\int_{x_0}^y \varphi'(y)^2e^{\varphi(y)}\p(X \ge y)dx dy \nonumber\\
&\le& 2\tilde{m}\E \int_{x_0}^{\tilde{m}}
\varphi'(y)^2e^{\varphi(y)}\ind_{\{X \ge y\}}dy \le 2\tilde{m}^2
\E\varphi'(X)^2 e^{\varphi(X)}. \label{positive_axes_3}
\end{eqnarray}
\paragraph{(ii)} $x_0 > \tilde{m}$: \quad We can obviously assume that
$\tilde{m} \ge 1$. Then as before
\begin{eqnarray}
&&
\int_{\tilde{m}}^{x_0}\int_0^{y}\varphi'(x)\varphi'(y)(e^{\varphi(x)}
+
e^{\varphi(y)})\p(X \le x)\p(X \ge y)dx dy \nonumber\\
&\le&
2\int_{\tilde{m}}^{x_0}\int_0^y\varphi'(x)^2e^{\varphi(x)}\p(X\le
x)\p(X
\ge y)dxdy \nonumber \\
&=& 2\int_{0}^{x_0}\int_{x\vee
\tilde{m}}^{x_0}\varphi'(x)^2e^{\varphi(x)}\p(X\le x)\p(X \ge y)dy
dx
\nonumber\\
&\le& 2\int_{0}^{x_0}\varphi'(x)^2e^{\varphi(x)}\p(X \le x)
\int_{\tilde{m}\vee x}^\infty y\p(X\ge
y)dy dx \nonumber\\
&\le& 2\sigma^2\int_{0}^{x_0}\varphi'(x)^2e^{\varphi(x)}\p(X \le
x)\p(X \ge x)dx \nonumber\\
&\le& 2\sigma^2\int_0^{\tilde{m}}\varphi'(x)^2e^{\varphi(x)}\p(X
\le x)dx +
2\sigma^2\int_{\tilde{m}}^{x_0}\varphi'(x)^2e^{\varphi(x)}x\p(X \ge x)dx \nonumber\\
&\le& 2\sigma^2(\tilde{m} +
2\sigma^2)\E\varphi'(X)^2e^{\varphi(X)}. \label{positive_axes_4}
\end{eqnarray}
Bringing together
(\ref{positive_axes_1}),(\ref{positive_axes_2}),(\ref{positive_axes_3})
and (\ref{positive_axes_4}) completes the proof. \qed

\begin{lemma}
\label{mixed_axes_lemma} Let $\varphi \colon \R \to \R$ be a
smooth convex Lipschitz function, non-increasing on $(-\infty,0)$,
and $X$ be a random variable with $\mathcal{L}(X) \in
\mathcal{M}(m,\sigma^2)$. Then
\begin{displaymath}
\int_{\{(x,y)\in \R^2 \colon xy\le
0\}}\varphi'(x)\varphi'(y)e^{\varphi(y)} \p(X \le x\wedge y)\p(X
\ge x\vee y)dxdy \le C(m,\sigma^2) \E\varphi'(X)^2e^{\varphi(X)}.
\end{displaymath}
\end{lemma}
\proof Assume without loss of generality that $m \ge 1$ and let
$\tilde{m}$ be the constant defined in Lemma
\ref{double_monotone_function}. For $x < 0 < y$ we have either
$\varphi'(x)\varphi'(y) \le 0$, or $\varphi'(x) \le \varphi'(y) <
0$ and $\varphi(x) \ge \varphi(y)$, so
\begin{eqnarray*}
&&\int_{-\infty}^0\int_0^\infty\varphi'(x)\varphi'(y)
e^{\varphi(y)}\p(X \le x)\p(X\ge y)dydx \\
&\le& \int_{-\infty}^0\int_0^\infty
\varphi'(x)^2e^{\varphi(x)}\p(X \le x)\p(X \ge y)dydx \le
C(m,\sigma^2) \int_{-\infty}^0\varphi'(x)^2e^{\varphi(x)}\p(X \le
x)dx,
\end{eqnarray*}
where the last inequality follows from the fact that by Lemma
\ref{subgaussian}
\begin{displaymath}
\int_0^\infty \p(X \ge y) dy= \E X_+ \le C(m,\sigma^2).
\end{displaymath}
Also
\begin{eqnarray*}
&&\int_0^\infty\int_{-\infty}^0
\varphi'(x)\varphi'(y)e^{\varphi(y)}\p(X \le y)\p(X \ge x) dy dx\\
&\le& \int_{-\infty}^0 \int_0^\infty
\varphi'(y)^2e^{\varphi(y)}\p(X \le y)\p(X \ge x)dx dy \le
C(m,\sigma^2)\int_{-\infty}^0 \varphi'(y)^2e^{\varphi(y)}\p(X \le
y)dy.
\end{eqnarray*}
Now
\begin{displaymath}
\int_{-\infty}^{-\tilde{m}}\varphi'(x)^2e^{\varphi(x)}\p(X \le
x)dx \le
\int_{-\infty}^{-\tilde{m}}\varphi'(x)^2e^{\varphi(x)}(-x)\p(X \le
x)dx \le 2\sigma^2\E\varphi'(X)^2e^{\varphi(X)}
\end{displaymath}
by Lemma \ref{double_monotone_function}, as $\mathcal{L}(X) \in
\mathcal{M}(m,\sigma^2)$. Moreover,
\begin{eqnarray*}
&&\int_{-\tilde{m}}^0 \varphi'(x)^2e^{\varphi(x)}\p(X \le x)dx =
\E \int_{-\tilde{m}}^0  \varphi'(x)^2e^{\varphi(x)}\ind_{\{X \le
x\}}dx \le \tilde{m}\E\varphi'(X)^2e^{\varphi(X)}.
\end{eqnarray*}
\qed

\paragraph{\noindent \textbf{Proof of Theorem \ref{entropy_estimates}\quad}}
We will follow Ledoux's approach for bounded variables. Due to the
tensorization property of entropy (Theorem
\ref{tensorization_theorem}), it is enough to prove the theorem
for $n=1$. Also, by the standard approximation argument, we can
restrict our attention to convex Lipschitz functions only. Let now
$Y$ be an independent copy of $X$. By Jensen's inequality we have
\begin{eqnarray*}
\Ent e^{\varphi(X)} &=& \E\varphi(X)e^{\varphi(X)} - \E
e^{\varphi(X)}\log \E e^{\varphi(X)} \le \frac{1}{2}\E(\varphi(X)
- \varphi(Y))(e^{\varphi(X)} - e^{\varphi(Y)}) \\
&=& \E(\varphi(X) - \varphi(Y))(e^{\varphi(X)} -
e^{\varphi(Y)})\ind_{\{X \le Y\}}\\
&=& \E \int_\R\int_\R \varphi'(x)\varphi'(y)e^{\varphi(y)}
\ind_{\{X \le x \le
Y\}}\ind_{\{X \le y \le Y\}} dx dy \\
& = & \int_\R\int_\R \varphi'(x)\varphi'(y)e^{\varphi(y)}\p(X \le
x\wedge y)\p(X \ge x\vee y)dx dy.
\end{eqnarray*}
Since $\mathcal{L}(-X)\in \mathcal{M}(m,\sigma^2)$, we can assume
that the minimal value of $\varphi$ is attained at some point of
the right semi-axis (possibly at $\infty$). Splitting now the
double integral into four integrals depending on the signs of $x$
and $y$ and using Lemmas \ref{positive_axes_lemma} and
\ref{mixed_axes_lemma} we obtain the desired inequality. Note that
we can use Lemma \ref{positive_axes_lemma} to handle the
integration over $(-\infty,0)^2$ again by change of variables and
the fact that $\mathcal{L}(-X) \in \mathcal{M}(m,\sigma^2)$. The
tail inequality follows from the entropy estimates by Theorem
\ref{log_sobolev_theorem}. \qed

\remark One would obviously like to characterize all real random
variables $X$ such that the random vectors $(X_1,\ldots,X_n)$
(where $X_i$'s are independent copies of $X$) have the
concentration property of order 2 for 1-Lipschitz convex functions
with constants independent of the dimension $n$. Each such
variable must of course have the concentration property of order 2
itself. This, however, is not sufficient, as concentration with
respect to convex functions implies hypercontractivity (see
\cite{Hitczenko}), which is equivalent to some regularity of the
tail. In particular it follows that $\E X^2\ind_{\{X \ge t\}} \le
Ct^2\p(X \ge t)$ for $t$ large. This condition is weaker than
$\mathcal{L}(X) \in \mathcal{M}(m,\sigma^2)$ for some $m,\sigma$
but hypercontractivity is also weaker than the concentration
property of order 2, uniformly over the dimension $n$.
\paragraph{}
We would also like  to point out that all Borel probability
measures $\mu$ on the real line, which satisfy the logarithmic
Sobolev inequality for all smooth (not necessarily log-convex)
functions, belong to $\mathcal{M}(m,\sigma^2)$ for some
$m,\sigma$. Thus $\cup_{m,\sigma}\mathcal{M}(m,\sigma^2)$ is
strictly larger than the class of all measures satisfying the
logarithmic Sobolev inequality for all smooth functions. More
precisely, we have the following
\begin{prop} Let $\mu$ be a Borel probability measure on $\R$ for
which there exists $C < \infty$ such that for all smooth
functions,
\begin{displaymath}
\Ent f(X)^2 \le C\E|f'(X)|^2,
\end{displaymath}
where $X$ is a random variable with the law $\mu$. Then there
exist constants $m, \sigma < \infty$ such that $\mu \in
\mathcal{M}(m,\sigma^2)$.
\end{prop}
\proof From the Bobkov-G\"otze criterion (see \cite{BG}) it
follows that if $n$ is the density of the absolutely continuous
part of $\mu$ and $M$ is a median of $\mu$, than for some constant
$K$ and all $x \ge M$,
\begin{displaymath}
\mu([x,\infty))\log\frac{1}{\mu([x,\infty))}\int_M^x
\frac{1}{n(t)}dt < K.
\end{displaymath}
Thus (since $\mu$ has the concentration property of order $2$)
from H\"older's inequality we get (using the above inequality for
$x + 1/x$ instead of $x$)
\begin{displaymath}
\frac{D}{K}x^2\mu([x+1/x,\infty))\le \frac{1}{\int_M^{x+1/x}
1/n(t)dt} \le \frac{\mu([x,x+1/x))}{(\int_x^{x+1/x} 1 dt)^2} =
x^2\mu([x,x+1/x)),
\end{displaymath}
for some $D > 0$ and $x$ large enough, which implies
$\mu([x+1/x,\infty)) \le \alpha\mu([x,\infty))$ with $\alpha =
\frac{K/D}{K/D + 1}$. Since a similar condition on the left tail
can be proven analogously, the claim follows by Proposition
\ref{equivalence}. \qed
\section{Concentration for seminorms on the tensor product and random chaoses \label{Sec_tensor_product}}
\subsection{A tensorization inequality for seminorms}
\begin{theorem}\label{seminorms_theorem}
Consider a random vector $X \in \R^n$ for which there exists a
constant $K$ such that for every seminorm $\varphi \colon \R^n \to
\R_+$ we have $\E \varphi(X) < \infty$ and for every $p \ge 1$,
\begin{displaymath}
\|\varphi(X) - \E \varphi(X)\|_p \le K\sqrt{p}\sup_{|\alpha| =
1}\varphi(\alpha).
\end{displaymath}
Then if $X_1,\ldots,X_d$ are independent copies of $X$ and $\psi
\colon \bigotimes_{i=1}^d \R^n \to \R_+$ a seminorm,  we have
\begin{displaymath}
\|\psi(\bigotimes_{i=1}^d X_i) - \E\psi(\bigotimes_{i=1}^dX_i)\|_p
\le K_d\sum_{I \subseteq\{1,\ldots,d\}, I \neq
\emptyset}p^{\#I/2}\E\sup_{|\alpha_k| \le 1 \colon k \in
I}\psi\left(\bigotimes _{i=1}^d X_{i,I,(\alpha_k)_{k\in
I}}\right),
\end{displaymath}
where
\begin{displaymath}
X_{i,I,(\alpha_k)_{k\in I}} = \left\{
\begin{array}{lcl}
X_i & {\rm if} & i \notin I, \\
\alpha_i & {\rm if} & i \in I,
\end{array}
\right.
\end{displaymath}
$K_d$ is a constant depending only on $K$ and $d$, and $|\alpha|$
denotes the Euclidean norm of the vector $\alpha$.
\end{theorem}

\proof We use induction on $d$. For $d=1$, the statement of the
theorem is just its hypothesis. Assume that the statement is true
for fewer than $d \ge 2$ copies of $X$. Using conditionally the
induction assumption for $d_1=1$ and the function
\begin{displaymath}
\varphi_1(x) = \psi(X_1\otimes\ldots\otimes X_{d-1}\otimes x),
\end{displaymath}
we obtain
\begin{displaymath}
\E_{X_d} \left| \psi(\bigotimes_{i=1}^d X_i) -
\E_{X_d}\psi(\bigotimes_{i=1}^d X_i) \right|^p \le K^p
p^{p/2}\sup_{|\alpha| \le 1} \psi(\bigotimes_{i=1}^{d-1}X_i
\otimes \alpha)^p.
\end{displaymath}
Now notice that the function $\varphi_2(x) = \sup_{|\alpha| \le 1}
\psi(x\otimes \alpha)$ is a seminorm on
$\bigotimes_{i=1}^{d-1}\R^d$ and thus we can apply the induction
assumption and the triangle inequality in $L^p$, which together
with the Fubini Theorem gives
\begin{equation} \label{tensor_chaos_1}
\E \left| \psi(\bigotimes_{i=1}^d X_i) -
\E_{X_d}\psi(\bigotimes_{i=1}^d X_i) \right|^p \le
\tilde{K}_{d-1}^p\sum_{I\subseteq\{1,\ldots,d\}, d \in I}
p^{p\#I/2}\left(\E\sup_{|\alpha_k| \le 1 \colon k \in I}
\psi(\bigotimes_{i=1}^d X_{i,I,(\alpha_k)_{k \in I}})\right)^p,
\end{equation}
where $\tilde{K}_{d-1}$ depends only on $K_{d-1}$ and $d$.

Now we would like to estimate $\E|\E_{X_d} \psi(\bigotimes_{i=1}^d
X_i) - \E\psi(\bigotimes_{i=1}^d X_i)|^p$. To this end consider
$\varphi_3 \colon\bigotimes_{i=1}^{d-1} \R^n\to \R_+$ defined as
$\varphi_3(x) = \E\psi(x\otimes X_d)$. It is easy to see that
$\varphi_3$ is a seminorm and thus, by the induction assumption,
\begin{equation} \label{tensor_chaos_2}
\E|\varphi_3(\bigotimes_{i=1}^{d-1} X_i) - \E
\varphi_3(\bigotimes_{i=1}^{d-1} X_i)|^p \le
\tilde{K}_{d-1}^p\sum_{I \subseteq \{1,\ldots, d-1\}, I \neq
\emptyset}p^{p\#I/2}\left(\E\sup_{|\alpha_k| \le 1 \colon k \in I}
\varphi_3\left(\bigotimes _{i=1}^{d-1} X_{i,I,(\alpha_k)_{k\in
I}}\right)\right)^p.
\end{equation}
Now it is enough to note that for each $I \subseteq
\{1,\ldots,d-1\}$,
\begin{displaymath}
\E\sup_{|\alpha_k| \le 1 \colon k \in I} \varphi_3\left(\bigotimes
_{i=1}^{d-1} X_{i,I,(\alpha_k)_{k\in I}}\right) \le
\E\sup_{|\alpha_k| \le 1 \colon k \in I} \psi\left(\bigotimes
_{i=1}^d X_{i,I,(\alpha_k)_{k\in I}}\right),
\end{displaymath}
which together with (\ref{tensor_chaos_1}) and
(\ref{tensor_chaos_2}) completes the proof. \qed

Notice that by Theorem \ref{entropy_estimates} for all product
random vectors $X \in \R^n$ with 1-dimensional marginals in
$\mathcal{M}(m,\sigma^2)$ and all seminorms $\varphi\colon \R^n
\to \R_+$ we have $\|(\varphi(X) - \E\varphi(X))_+\|_p \le K
\sqrt{p}\sup_{|\alpha| \le 1}\varphi(\alpha)$, with $K$ depending
only on $m$ and $\sigma^2$. Thus the same proof, with formal
changes only, gives
\begin{theorem}\label{seminorms_theorem_upper_tail}
Let $X_1,\ldots,X_n \in \R^n$ be independent random vectors with
independent components and all 1-dimensional marginals in
$\mathcal{M}(m,\sigma^2)$. Then for every seminorm $\psi \colon
\bigotimes_{i=1}^d \R^n \to \R_+$ we have
\begin{displaymath}
\|(\psi(\bigotimes_{i=1}^d X_i) - \E\psi(\bigotimes_{i=1}^d
X_i))_+\|_p \le K_d\sum_{I \subseteq\{1,\ldots,d\}, I \neq
\emptyset}p^{\#I/2}\E\sup_{|\alpha_k| \le 1 \colon k \in
I}\psi\left(\bigotimes _{i=1}^d X_{i,I,(\alpha_k)_{k\in
I}}\right),
\end{displaymath}
where $K_d$ depends only on $m,\sigma^2$ and $d$.
\end{theorem}

By the Chebyshev inequality we can obtain from the above theorems
a corollary concerning the tail behaviour of
$\psi(\bigotimes_{i=1}^d X_i)$. We give it only for Theorem
\ref{seminorms_theorem}; for Theorem
\ref{seminorms_theorem_upper_tail} it is analogous but deals with
the upper tail only.
\begin{cor}
Under the assumption of Theorem \ref{seminorms_theorem} there
exist constants $K_d$, depending only on $K$ and $d$, such that
for all $t \ge 1$,
\begin{displaymath}
\p\left( |\psi(\bigotimes_{i=1}^d X_i) - \E\psi(\bigotimes_{i=1}^d
X_i)|\ge K_d\sum_{I \subseteq\{1,\ldots,d\}, I \neq
\emptyset}t^{\#I/2}\E\sup_{|\alpha_k| \le 1 \colon k \in
I}\psi\left(\bigotimes _{i=1}^d X_{i,I,(\alpha_k)_{k\in
I}}\right)\right) \le e^{-t}.
\end{displaymath}
\end{cor}

\subsection{Chaos random variables} The above theorems can be
rewritten in terms of decoupled polynomial chaoses. Below we
present such a version as a corollary (actually equivalent to
Theorem \ref{seminorms_theorem}), so that the reader could compare
it with existing results on chaos random variables.

Let $X^{(1)},\ldots,X^{(d)}$ be independent copies of
$X=(X_1,\ldots,X_n)\in \R^n$ and consider a homogeneous decoupled
chaos of order $d$, i.e. a random variable of the form
\begin{equation}
Z  = \sup_{t \in \mathcal{T}} \left|\sum_{i_1,\ldots,i_d=1}^n
t_{i_1\ldots i_d} X_{i_1}^{(1)}\ldots X_{i_d}^{(d)}\right|,
\label{chaos_definition}
\end{equation}
where $\mathcal{T}$ is a countable, bounded set of functions $t
\colon \{1,\ldots,n\}^d \to \R$.

Let us introduce
\begin{defi}
\label{coefficients_definition} For $I \subseteq \{1,\ldots,d\}$
let
\begin{displaymath}
\|\mathcal{T}\|_I = \E \sup_{t \in \mathcal{T}}
\sup_{|\alpha^{(k)}| \le 1, k \in I}\left| \sum_{i_1\ldots i_d =
1}^n t_{i_1,\ldots, i_d} \prod_{k \in I}\alpha_{i_k}^{(k)}\prod_{k
\notin I}X_{i_k}^{(k)}\right|,
\end{displaymath}
where the second supremum is taken over all $(\alpha^{(k)})_{k\in
I} \in (\R^n)^{\#I}$, and $|\alpha^{(k)}|$ stands for the
Euclidean norm of the vector $\alpha^{(k)}$.
\end{defi}

\begin{cor}
\label{Cor_chaoses} Let $X = (X_1,\ldots,X_n)$ have the
concentration property with constants $C,K$ with respect to
1-Lipschitz seminorms on $\R^n$. Then for any integer $d \ge 1$,
there exists a constant $K_d$, depending only on $d$ and $C,K$,
such that for any homogeneous chaos Z (as defined in
(\ref{chaos_definition})) and any $p \ge 1$,
\begin{equation}
\|Z - \E Z\|_p \le K_d\sum_{I \subseteq \{1,\ldots,d\}, I \neq
\emptyset} p^{\#I/2}\|\mathcal{T}\|_I.
\end{equation}
\end{cor}

\begin{cor}
There exist constants $K_d$ such that for all $t\ge 1$,
\begin{displaymath}
\p\left(|Z - \E Z| \ge K_d \sum_{I\subseteq \{1,\ldots,d\}, I\neq
\emptyset} t^{\# I/2}\|\mathcal{T}\|_I\right) \le e^{-t}.
\end{displaymath}
\end{cor}

\remark Usually one is interested in the 'undecoupled' chaos, i.e.
a random variable of the form
\begin{displaymath}
Z = \sup_{t\in \mathcal{T}} \left|\sum_{i_1,\ldots,i_d=1}^n
t_{i_1\ldots i_d}X_{i_1}\ldots X_{i_d}\right|,
\end{displaymath}
where $X_1,\ldots,X_n$ are independent random variables and for
all $t\in \mathcal{T}$ the number $t_{i_1\ldots i_d}$ is invariant
under permutations of coordinates and non-zero only if the
coordinates are pairwise distinct. The analogues of the above
theorems for such chaoses generated by  Gaussian variables were
obtained by Borell \cite{BO} and Arcones, Gin\'{e} \cite{Arcones}.
The Rademacher case was considered by Talagrand \cite{T2} (chaos
of order 2) and Boucheron, Bousquet, Lugosi and Massart
\cite{BBLM} (chaos of arbitrary order $d$). It is easy to see that
each decoupled chaos can be represented as an undecoupled one, but
the aforementioned results do not recover the lower tail
inequalities (except for the case $d=2$). Also the methods are
quite different and do not allow treating both Gaussian and
Rademacher variables in a unified way.

\subsubsection{Chaoses generated by symmetric random variables
with log-concave tails}

Now we would like to point out that the proof of Theorem
\ref{seminorms_theorem} can be actually used in a slightly
different setting, namely for chaoses generated by independent
random variables with logarithmically concave tails. Such
variables have been investigated by Lata\l a \cite{Latala3} and \L
ochowski \cite{Lochowski}.
\begin{defi} Let $\mathcal{N} = (X_i^{(k)})_{k\le d, i\le n}$ be a
matrix of independent symmetric random variables with
logarithmically concave tails, i.e. random variables such that the
functions
\begin{displaymath}
\mathcal{N}_i^{(k)}(t) = -\log \p(|X_i^{(k)}| \ge t), \;t \ge 0
\end{displaymath}
are convex. Furthermore assume (as a matter of normalization) that
\begin{displaymath}
\inf\{t\colon \mathcal{N}_i^{(k)}(t) \ge 1\} = 1,
\end{displaymath}
and define modified functions $\tilde{\mathcal{N}}_i^{(k)}$ by the
formula
\begin{displaymath}
\tilde{\mathcal{N}}_i^{(k)}(t) = \left\{
\begin{array}{ccl}
t^2& {\rm for} & |t| \le 1\\
\mathcal{N}_i^{(k)}(|t|) & {\rm for} & |t| \ge 1.
\end{array}\right.
\end{displaymath}
\end{defi}

Let now $\mathcal{T}$ be a countable set of functions $t\colon
\{1,\ldots,n\}^d \to \R$ and $Z$ a random variable defined by
(\ref{chaos_definition}). Moment estimates for $Z$ will be
expressed in terms of the following quantities:

\begin{defi} For $I \subseteq \{1,\ldots,d\}$ and $p \ge 1$ define
\begin{displaymath}
\|\mathcal{T}\|_{\mathcal{N},I,p} = \E \sup_{t \in
\mathcal{T}}\;\sup_{\alpha^{(k)}\in \mathcal{A}_{k,p} \colon k \in
I}\; |\sum_{i_1,\ldots,i_d = 1}^nt_{i_1\ldots i_d}\prod_{k \in
I}\alpha_{i_k}^{(k)}\prod_{k \notin I}X_{i_k}^{(k)}|,
\end{displaymath}
where
\begin{displaymath}
\mathcal{A}_{k,p} = \{\alpha \in \R^n \colon \sum_{i=1}^n
\tilde{\mathcal{N}}_i^{(k)}(\alpha_i) \le p\}.
\end{displaymath}
\end{defi}

The following result was proved for $d=1$ by Lata\l a
\cite{Latala3} and for arbitrary $d$ by \L ochowski
\cite{Lochowski}.
\begin{theorem}[Lata\l a, \L ochowski] \label{Latala_Lochowski}There exist constants
$K_d$ (depending only on $d$) such that for all $p \ge 1$,
\begin{displaymath}
\frac{1}{K_d}\sum_{I\subseteq\{1,\ldots,d\}}\|\mathcal{T}\|_{\mathcal{N},I,p}\le\|Z\|_p
\le K_d\sum_{I\subseteq
\{1,\ldots,d\}}\|\mathcal{T}\|_{\mathcal{N},I,p}.
\end{displaymath}
\end{theorem}

\begin{prop}
The conclusion of Theorem \ref{Latala_Lochowski} for $d=1$ implies
this conclusion for arbitrary $d$.
\end{prop}
\proof When we rewrite the inequalities of Theorem
\ref{Latala_Lochowski} in the language of seminorms, the proof
becomes analogous to the proof of Theorem \ref{seminorms_theorem}.
One has simply to notice that the factors $p^{\#I/2}$ do not
appear, as the dependence on $p$ is incorporated in the sets
$\mathcal{A}_{k,p}$ which replace the unit Euclidean ball in the
supremum.

 As for the lower estimate, the proof is even simpler but we
present it here for the sake of completeness (written in the
'chaos language'). Obviously $\|Z\|_p \ge \|Z\|_1=
\|\mathcal{T}\|_{\mathcal{N},\emptyset,p}$. Moreover for any
nonempty set $I\subseteq\{1,\ldots,d\}$, say $r \in I$, we have by
the induction hypothesis
\begin{displaymath}
\|Z\|_p \ge \frac{1}{K_1}\left\|\sup_{t \in
\mathcal{T}}\sup_{\alpha^{(r)}\in
\mathcal{A}_{r,p}}\left|\sum_{i_1,\ldots,i_d=1}^n t_{i_1\ldots
i_d}\prod_{k\neq r}X_{i_k}^{(k)}\alpha_{i_r}^{(r)} \right|
\right\|_p \ge \frac{1}{K_1
K_{d-1}}\|\mathcal{T}\|_{\mathcal{N},I,p}.
\end{displaymath}
 \qed

Using the Chebyshev inequality and the Paley-Zygmund inequality
together with the hypercontractive properties of chaoses (see
\cite{dlp2}) we obtain
\begin{cor}
\label{tail_estimates_log_concave} There exist constants $K_d$
such that for all $t \ge 1$,
\begin{eqnarray*}
\p\left(Z \ge K_d \sum_{I \subseteq
\{1,\ldots,d\}}\|\mathcal{T}\|_{\mathcal{N},I,t}\right) &\le& e^{-t} \\
\p\left(Z \ge \frac{1}{K_d} \sum_{I \subseteq
\{1,\ldots,d\}}\|\mathcal{T}\|_{\mathcal{N},I,t}\right) &\ge&
e^{-t}\wedge c.
\end{eqnarray*}
\end{cor}
\remark  The bound of Theorem \ref{Latala_Lochowski} is also valid
in the case of undecoupled chaoses due to the decoupling results
by de la Pe\~{n}a and Montgomery-Smith \cite{dlp1}, which say that
tails and moments of the decoupled and undecoupled chaos with the
same coefficients are equivalent.

\remark All estimates presented so far are expressed in terms of
expected values of empirical processes, which themselves are in
general troublesome and difficult to estimate. One would obviously
want to obtain moment estimates in terms of deterministic
quantities at least in the real-valued case (i.e. when
$\mathcal{T}$ is a singleton). It has been done by Lata\l a in
\cite{Latala3} for $d=2$ and log-concave random variables and
recently in \cite{Latala5} for arbitrary $d$ and Gaussian chaoses.

\section{An application \label{Sec_application}}
Finally, we would like to argue that estimates in the spirit of
Section \ref{Sec_tensor_product}, although non-deterministic, may
be of some use. We will demonstrate it by presenting a sketch of a
new (at least to the author's best knowledge) proof of exponential
integrability of generalized Rademacher chaos process for general
$d$, which we believe is simpler than the preceding ones. The
general result and the proof for $d \le 2$ may be found in the
monograph by Ledoux and Talagrand \cite{LT}.

Let us first recall the general setting. We deal with a Banach
space $B$ for which there exists a countable set $D$ of linear
functionals from the unit ball of $B'$ such that for each $x \in
B$ we have $\|x\| = \sup_{f \in D} |f(x)|$. A $B$-valued random
variable $X$ is a homogeneous Rademacher chaos of order $d$ if
there is a sequence $(x_{i_1\ldots i_d})_{i_1,\ldots,i_d \in \N}
\in B$ ($x_{i_1\ldots i_d}$ invariant under permutation of
coordinates and non-zero only if $i_1,\ldots,i_d$ are pairwise
distinct) such that for every $f \in D$ the multiple series
$\sum_{i_1,\ldots,i_d}f(x_{i_1\ldots
i_d})\varepsilon_{i_1}\ldots\varepsilon_{i_d}$ converges almost
surely and $(\sum_{i_1,\ldots,i_d}f(x_{i_1\ldots
i_d})\varepsilon_{i_1}\ldots\varepsilon_{i_d})_{f\in D}$ has the
same distribution as $(f(X))_{f\in D}$.

\begin{theorem}\label{exp_integ}
Let $X$ be a homogeneous Rademacher chaos of order $d$. Then for
all $\alpha \in \R$
\begin{displaymath}
\E e^{\alpha\|X\|^{2/d}} < \infty.
\end{displaymath}
\end{theorem}

Before we proceed with the proof, we need
\begin{lemma}
\label{Rademacher_chaos_lemma} There exist constants $L_d$
(depending only on $d$) with the property that for every $\alpha >
0$ there exists a constant $\varepsilon(\alpha)$ such that for
every homogeneous Rademacher chaos $X$ of order $d$ and every $M$
satisfying $\p(\|X\|> M) < \varepsilon(\alpha)$ we have
\begin{displaymath}
\p(\|X\| > L_dMt^{d/2}) < e^{-\alpha t}
\end{displaymath}
for all $t \ge 1$.
\end{lemma}

\proof Rademacher variables are log-concave, so we can apply
Corollary \ref{tail_estimates_log_concave} for finite sums. It is
easy to see that
\begin{equation}\label{exp_integ_lemma_sets}
\mathcal{A}_{k,p} = \mathcal{A}_p = \{(\alpha_i)\colon
\sum\alpha_i^2 \le p, |\alpha_i| \le 1\}.
\end{equation}
Due to convergence, using standard arguments one can extend the
tail estimates to the general Rademacher chaos $X$. Therefore, set
\begin{displaymath}
\phi(t) =
\sum_{I\subseteq\{1,\ldots,d\}}\|\mathcal{T}\|_{\mathcal{N},I,t},
\end{displaymath}
where $ \mathcal{T} = \{(f(x_{i_1\ldots i_d}))\colon f \in D\}$.
From (\ref{exp_integ_lemma_sets}) it follows that
\begin{equation}\label{exp_integ_lemma_ineq}
\phi(xt) \le t^{d/2}\phi(x)
\end{equation}
for all $x \ge 0, t \ge 1$. As $\phi$ is increasing, by the last
inequality it is continuous. If $\phi(x) \le K_d M$ for every $x$
($K_d$ being the constant from Corollary
\ref{tail_estimates_log_concave}), we have, for $t \ge 1$,
\begin{displaymath}
\p(\|X\| \ge MK_d^2t^{d/2}) \le \inf_x \p(\|X\| \ge K_d \phi(x))
\le \inf_x e^{-x} = 0.
\end{displaymath}
Otherwise $K_d M = \phi(x)$ for some $x$. Thus
\begin{displaymath}
\varepsilon \ge \p(\|X\| \ge M) = \p(\|X\| \ge \phi(x)/K_d) \ge
c\wedge e^{-x},
\end{displaymath}
which for $\varepsilon$ small enough yields $x \ge -\log
\varepsilon$. Moreover (\ref{exp_integ_lemma_ineq}) gives
$K_dMt^{d/2} \ge \phi(tx)$ for $t \ge 1$ and thus for $\varepsilon
< e^{-\alpha}$,
\begin{displaymath}
\p(\|X\| \ge K_d^2 M t^{d/2}) \le \p(\|X\| \ge K_d\phi(tx)) \le
e^{-tx} \le e^{t\log\varepsilon} \le e^{-\alpha t}.
\end{displaymath}
\qed

\noindent{\textbf{Proof of Theorem \ref{exp_integ}} (sketch)\quad}
We will proceed by induction on $d$.  Let $S_n = \sum_{i_1
,\ldots,i_d \ge n} x_{i_1\ldots
i_d}\varepsilon_{i_1}\ldots\varepsilon_{i_d}$. Then $\|S_n\|$ is a
reversed submartingale with $\E\| S_n\| \le \E\|X\|$ and thus it
converges to some random variable, which by the zero-one law must
be almost surely constant. Thus there exists $M$ such that for all
$\varepsilon > 0$ we have $\p(\|S_n\| \ge M) < \varepsilon$ for
$n$ large enough. By Lemma \ref{Rademacher_chaos_lemma}, for every
$\beta > 0$ there is $n$ such that for all $t \ge 1$,
\begin{displaymath}
\p(\|S_n\| \ge L_dM t^{d/2}) < e^{-\beta t}
\end{displaymath}
or equivalently
\begin{displaymath}
\p(\alpha\|S_n\|^{2/d} \ge t) \le \exp\left(-\frac{\beta t
}{\alpha(L_d M)^{2/d}}\right)
\end{displaymath}
for $t \ge (L_d M)^{2/d}\alpha$, which clearly implies $\E
e^{\alpha\|S_n\|^{2/d}} < \infty$ for $\alpha < \beta/(L_d
M)^{2/d}$. Since $X - S_n$ is a finite sum of chaoses of orders
lower than $d$ (when $d > 1$) or a bounded random variable (for $d
= 1$), its integrability properties allow us to use H\"older's
inequality and obtain $\E e^{\alpha \|X\|^{2/d}} < \infty$ also
for $\alpha < \beta/(L_dM)^{d/2}$. This allows us to finish the
proof as $\beta$ can be chosen arbitrarily large. \qed

\paragraph{Acknowledgements} The author would like to express his
gratitude to Prof. Rafa\l\;Lata\l a for introducing him to the
subject and all the inspiring conversations.

Rados\l aw Adamczak \\
Institute of Mathematics \\
Polish Academy of Sciences \\
 00-956 Warszawa 10, Poland\\
ul. Śniadeckich 8. P.O.Box 21\\
e-mail: R.Adamczak@impan.gov.pl
\end{document}